\documentclass[11pt]{article}
\usepackage{amsfonts,amssymb}
\usepackage{amsmath}
\usepackage{amsthm}
\usepackage[matrix,arrow,curve]{xy}
\newtheorem{theorem}{Theorem}
\newtheorem{lemma}{Lemma}

\newtheorem{corollary}{Corollary}

\title{Soft Construction of Floer-type Homologies}
\author{A.~A.~Agrachev\thanks{SISSA, Trieste and Steklov. Math. Inst., Moscow}}
\date{}

\begin{document}
\maketitle

\begin{abstract}
Floer homology is a good example of homological invariants living in the infinite dimension. We suggest a way to construct
this kind of invariants using only soft essentially finite-dimensional tools; no hard analysis or PDE is involved.
This work is partially inspired by the M. Gromov's survey ``Soft and hard symplectic geometry'' at ICM-86.
\end{abstract}

\section{Introduction}

Floer homology (see \cite{Fl,Fl1}) is a good example of homological invariant, which lives in the infinite dimension. It is not completely homotopy invariant\footnote{indeed, infinite dimension implies an infinite freedom :-)}, but it survives a class of ``compact homotopies'' that are close enough to finite-dimensional perturbations.

Classical construction of Floer homology and similar invariants involves hard analysis and elliptic systems of PDEs.
This is a beautiful theory but it is technically heavy and rather restrictive; we cannot go beyond a limited number of situations where good elliptic systems are available.

In this expository paper, partially motivated by the well-known M.~Gromov's survey \cite{Gr}, we suggest a soft and essentially
finite-dimensional approach to the Floer-type homological invariants. We use finite-dimensional subspaces of an auxiliary Banach space for finite-dimensional approximations of the
original big object. The role of this Banach space is similar to the role
of the events space in the Kolmogorov's foundation of Probabilities Theory.

In our study, the infinite dimension and infinite size remain only potential infinities. We deal with very big
but still finite-dimensional and compact objects, the dimension and size are just big parameters. The ``Infinity'' is
treated as a singularity and we study asymptotics of usual finite-dimensional homological invariants at this singularity.

Let me describe a toy example, the Leray--Schauder degree. Let $B$ be an infinite-dimensional separable Banach space and
$S\subset B$ the unite sphere in $B$. Let $\mathcal E=\{E\subset B: \dim E<\infty\}$ be the ordered by the inclusion directed set of finite-dimensional subspaces of $B$. Given $E\in\mathcal E$, the intersection $E\cap S$ is a
$(\dim E-1)$-dimensional sphere. In particular, the homology group $H_{i(\dim E-1)}(S\cap E)$ equals $\mathbb Z$ for $i=0,1$ and
equals 0 for other $i$; it depends only on $i$ and does not depend on $E$. We set:
$$
G_i(S)=H_{i(\dim E-1)}(S\cap E)
$$
and call $G_i(S)$ the Leray--Schauder homology of $S$. Recall that $S$ is contractible and $H_i(S)=0,\ \forall i\ne 0$.

Now let $\varphi:S\to B$ be a compact map such that $x+\varphi(x)\ne 0,\ \forall x\in S$. For any $\varepsilon>0$, there
exist a $\varepsilon$-close to $\varphi$ finite-dimensional map $\varphi_\varepsilon:S\to E_\varepsilon$, where
$E_\varepsilon\in\mathcal E$. If $\varepsilon$ is small enough, then $x+\varphi_\varepsilon(x)\ne 0,\ \forall x\in S$, that
allows us to define a map:
$$
\Phi^\varepsilon_E:S\cap E\to S\cap E,\quad \Phi^\varepsilon_E(x)=\frac{x+\varphi_\varepsilon(x)}{|x+\varphi_\varepsilon(x)|},
$$
for any $E\supset E_\varepsilon$. The degree of this map $d=\deg(\Phi^\varepsilon_E)$ does not depend on $E$ and is the same for all sufficiently good approximations $\varphi_\varepsilon$. This is the Leray--Schauder degree.

The degree is defined by the homomorphism
$$
{\Phi^\varepsilon_E}_*:H_{\dim E-1}(S\cap E)\to H_{\dim E-1}(S\cap E), \quad {\Phi^\varepsilon_E}_*(c)=cd,
$$
for any $c\in H_{\dim E-1}(S\cap E)=\mathbb Z$. Since this homomorphism does not depend on $\varepsilon$ and $E$ (for small enough $\varepsilon$ and big enough $E$), we may interpret it as a homomorphism
$\Phi_*:G_1(S)\to G_1(S)$, where $\Phi=\frac{I+\varphi}{|I+\varphi|}$. This homomorphism survives homotopies $\varphi_t,\
0\le t\le 1,$ where all $\varphi_t$ are compact and $x+\varphi_t(x)\ne 0,\ \forall x\in S$.

In the next section we give our soft construction of Floer homology and in Section~3 we explain, why it works.
In Section~4 we briefly study another example, inspired by the sub-Riemannian geometry. This example demonstrates a much more interesting asymptotic behavior of the homology than in the usual Floer case. We conclude with a draft of general
construction that can be adapted to many different situations.

In what follows, the definitions and statements are precise while the proofs are only sketched.

\medskip {\sl Acknowledgments.} I am greatful to Antonio Lerario and Michele Stecconi for very stimulating discussions.
In particular, Michele realized that original version of Theorem~1 with 4 consecutive limits can be simplified and 3 limits are sufficient. I hope that Antonio, Michele and maybe other interested people will develop the designated here method.

\section{Floer Homology}

We consider a compact smooth manifold $M$ endowed with a symplectic structure $\sigma$. Let $\tilde M$ be the universal covering of $M$ and $\tilde\sigma$ the pullback of $\sigma$ to $\tilde M$; we assume that $\tilde\sigma$ is an exact form: $\tilde\sigma=ds$.

We denote by $\Omega$ the space of contractible closed curves in $M$ of class $H^1$. In other words, $\Omega$ consists of contractible maps $\gamma:S^1\to M$, where $\gamma$ is differentiable almost everywhere with the derivative of class $L^2$. The lifts of $\gamma\in\Omega$ to $\tilde M$ are closed curves and we use the same symbol $\gamma$ for any lift of this curve to $\tilde M$.

Let $h_t:M\to\mathbb R,\ t\in S^1,$ be a measurable bounded with respect to $t\in S^1$ family of smooth functions on $M$.
We are going to study the functionals $\varphi_h:\Omega\to\mathbb R$ defined by the formula
$$
\varphi_h(\gamma)=\int_{S^1}s(\dot\gamma(t))-h_t(\gamma(t))\,dt.
$$
Note that $\int_{S^1}\langle s,\dot\gamma(t)\rangle\,dt=\int_\Gamma\sigma$ for any film $\Gamma$ such that $\gamma=\partial\Gamma$ and this integral does not depend on the choice of the lift of $\gamma$ to $\tilde M$.
Given $c\in\mathbb R$, we denote by $\Omega^c_h$ the Lebesgue set of $\varphi_h$:
$$
\Omega^c_h=\{\gamma\in\Omega:\varphi_h(\gamma)\le c\}
$$

We use some auxiliary objects. The homology invariants will be constructed with a help of these objects but do not depend on their choice. Namely, we assume that $M$ is equipped with a Riemannian structure $\langle\cdot,\cdot\rangle$ adapted to the symplectic structure, i.\,e. $\sigma(\xi,\eta)=\langle J\xi,\eta\rangle,\ \xi,\eta\in TM,$ where $J:TM\to TM$ is a quasi-complex structure, $J^2=-I$.

Moreover, we fix generators $X_1,\ldots,X_l$ of the $C^\infty(M)$-module $\mathrm{Vec}M$ of all smooth vector fields on $M$
and define a linear map $X_q:\mathbb C^l\to T_qM$ by the formula
$$
X_qu=\sum_{j=1}^lv^jX_j(q)+w^jJX_j(q),
$$ where $u=(u^1,\ldots,u^l),\ u_j=v_j+iw_j\in\mathbb C,\ j=1,\ldots,l$. We assume that $\langle\cdot,\cdot\rangle|_{T_qM},\ q\in M,$ is the image of the standard Euclidean structure on $\mathbb C^l$ by the linear map $X_q$; in other words,
$$
\langle\xi,\xi\rangle=\min\{|u|^2 : u\in\mathbb C^l,\ \xi=X_qu\}.
$$

A simple model example is an even-dimensional torus with a constant sympletic structure, flat Riemannian metric and constant vector fields \linebreak $X_1,\ldots, X_l$. Needless to say that, in general, if $M$ is not parallelizable, number $l$ is greater than $\dim M$. Anyway, the described auxiliary objects always exist and can be easily built.

Let $W$ be the space of all curves in $M$ of class $H^1$ parameterized by the segment $[0,1]$. We fix a parametrisation of $S^1$ by $[0,1]$; then $\Omega\subset W$.

We define the map $\phi:M\times L^2([0,1];\mathbb C^l)\to W$ as follows. Given $q\in M$ and
$u(\cdot)\in L^2([0,1];\mathbb C^l)$ the curve $\gamma(\cdot)=\phi(q,u(\cdot))$ is the solution of the ordinary differential equation
$$
\dot\gamma(t)=X_{\gamma(t)}u(t), \quad 0\le t\le 1,
$$
with the initial condition $\gamma(0)=q$. It is important for us that weak convergence of a sequence $u_n\in L^2([0,1];\mathbb C^l)$ to $u$ ($u_n\rightharpoonup u$ as $n\to\infty$) and convergence of $q_n\in M$ to $q$ imply uniform convergence of $\phi(q_n,u_n)$ to $\phi(q,u)$.

We also set $\phi_t(q,u)=(q,\phi(q,u)(t))$ and thus define the map \linebreak $\phi_t:M\times L^2([0,1];\mathbb R^l)\to M\times M$.
It is easy to see that $\phi_t$ is a smooth map and $\phi_t$ is a submersion for $0<t\le 1$.

The differential of $\phi_t$ is a linear operator $D_{(q,u)}\phi_t:T_qM\times L^2([0,1];\mathbb C^l)\to T_qM\times T_qM$.
It is important for us that convergence $q_n\to q$ and weak convergence $u_n\rightharpoonup u$ as $n\to\infty$ imply convergence $D_{(q_n,u_n)}\phi_t \to D_{(q,u)}\phi_t$ in the operator norm. It follows immediately from the standard ODE's ``variations formula'' for the differental of $\phi_t$ and the uniform convergence $\phi(q_n,u_n)\to\phi(q,u)$.

Given a subspace $E\subset L^2([0,1];\mathbb C)$, we set $E^l=E\times\cdots\times E\subset L^2([0,1];\mathbb C^l)$.

\begin{lemma} $0<t\le 1,r>0$, and $B_r=\left\{u\in L^2([0,1];\mathbb C^l):\|u\|<r\right\}$. Then there exists a finite-dimensional subspace $E_r\subset L^2([0,1];\mathbb C)$ such that for any subspace $E\supset E_r$ the map
$\phi_t\bigr|_{M\times(E^l\cap B_r)}$ has no critical points.
\end{lemma}
{\bf Proof.} For any $(q,u)\in M\times L^2([0,1];\mathbb C^l)$ there exists a finite-dimensional subspace $E_{(q,u)}$
such that the restriction of $D_{(q,u)}\phi_t$ to the subspace \linebreak $T_qM\times E^l(q,u)$ is surjective; simply because $\phi_t$ is a submersion and $M$ is finite-dimensional.

Moreover, we can take $E_{(q',u')}=E_{(q,u)}$ for all $q'$ close to $q$ and $u'$ close to $u$ in the weak topology. Hence, due to the compactness of $M$ and weak compactness of $\bar B_r$, we may assume that
$\#\{E_{(q,u)}: q\in M,\ u\in B_r\}<\infty.$
Then $E_r=\sum\limits_{(q,u)\in M\times B_r}E_{(q,u)}$ is the desired space.\qquad $\square$

If $E$ satisfies conditions of Lemma~1, then
$$
U_r(E)\doteq\{(q,u)\in M\times (B_r\cap E^l): \phi_t(q,u)=(q,q)\}  \eqno (1)
$$
is a smooth  submanifold of $M\times E^l$ of dimension $\dim E^l$ and $\phi(U_r(E^l))\subset\Omega$.

In the next construction, we prefer to treat $S^1$ in the intrinsic way as an Abelian Lie group equipped with the Haar measure $d\theta$ rather than directly parameterize it by a segment. Let $q\in M$, we denote by $b_q$ a symmetric bilinear form on the space $H^1(S^1;T_qM)$ defined by the formula:
$$
b_q(\xi,\eta)=\int_{S^1}\sigma(\xi(\theta),\dot\eta(\theta))\,d\theta,\quad \xi,\eta\in H^1(S^1;T_qM).
$$
A vector-function $\eta\in H^1(S^1;T_qM)$ belongs to the kernel of $b_q$ if and only if $\eta$ is constant, i.\,e.
$\eta(\theta)\equiv\eta(0)$. Indeed, if $\dot\eta\ne 0$, then
$$
b_q(J\dot\eta,\eta)=-\int_{S^1}\langle\dot\eta(\theta),\dot\eta(\theta)\rangle\,d\theta<0. \eqno (2)
$$
So, if $\dot\eta\in H^1(S^1;T_qM)$, then $\eta$ is not in the kernel; otherwise, we approximate $\dot\eta$ by a smooth function in the norm $L^2$ and plug-in the approximating function in (2) instead of $\dot\eta$.

We denote by $\imath:H^1(S^1;T_qM)\to H^1(S^1;T_qM)$ the involution defined by the formula $(\imath\xi)(\theta)=\xi(-\theta)$.
Then
$$
b_q(\imath\xi,\imath\eta)=-b_q(\xi,\eta),\quad \xi,\eta\in H^1(S^1;T_qM).  \eqno (3)
$$

Let $E$ be a finite-dimensional subspace of $L^2([0,1];\mathbb C)$ and \linebreak $E_0=\left\{\upsilon\in E: \int_0^1\upsilon(t)\,dt=0\right\}$.
We set:
$$
\mathcal X_q(E)=\left\{\theta:\mapsto\xi_0+\int_0^\theta X_qu(t)\,dt: \xi_0\in T_qM,\ u(\cdot)\in E^l_0\right\}\subset H^1(S^1;T_qM).
$$

We say that $E$ is {\it well-balanced} if $\imath E=E$ and $\ker b_q|_{\mathcal X_q(E)}=\ker b_q$.

\begin{lemma} Any finite-dimensional subspace of $L^2([0,1];\mathbb C)$ is contained in a well-balanced subspace.
\end{lemma}
Indeed, we can always add to $E$ a big enough finite-dimensional space of trigonometric polynomials to guarantee that for the enlarged space $\hat E$  and any nonconstant $\eta\in\mathcal X_q(\hat E)$ there exists $\xi\in\mathcal X_q(\hat E)$ that is not orthogonal to $\dot\eta$ in $L^2(S^1;T_qM)$.

Let $\mathcal E$ be the directed set of well-balanced subspaces ordered by the inclusion, the $\mathcal E\mbox{-}\lim$ of a generalized sequence indexed by the elements of $\mathcal E$ is defined in a usual way.

We are now ready to formulate the main result. Let $E\in\mathcal E,\ c>0$, and $r>0$. We consider relative homology groups
$$
G_i(E;c,r)\doteq H_i\left(\phi(U_r(E))\cap\Omega^c_h,\phi(U_r(E))\cap\Omega^{-c}_h\right),
$$
$i=0,1,2,\ldots$.

Let $\jmath_r:G_i(E;c,r)\to H_i\left(\phi(E)\cap\Omega^c_h,\phi(E)\cap\Omega^{-c}_h\right)$ be the homology homomorphisms induced by the inclusions
$$
\phi(U_r(E))\cap\Omega^{\pm c}_h\subset\phi(E)\cap\Omega^{\pm c}_h.
$$

\begin{theorem}
There exist
$$
\lim\limits_{c\to\infty}\lim\limits_{r\to\infty}\mathcal E\mbox{-}\lim \jmath_r\bigl(G_{i+d_E}(E;c,r)\bigr)\cong H_i(M),
$$
where $d_E=\frac 12(\dim E-1)\dim M$.
\end{theorem}

{\bf Remark.} Actually, the families under the limits stabilize; a little bit more precise statement is as follows.
For any big enough $c>0$ there exist $r(c)>0$ and $E(c,r)\in\mathcal E$ such that for any $r>r(c)$ and well-balanced $E\supset E(c,r)$ we have:
$$
\jmath_r\left(G_{i+d_E}(E;c,r)\right)\cong H_i(M).
$$

The limit in Theorem~1 is a soft construction of Floer homology. It implies Morse inequalities for the functional $\varphi_h$, which we are going to describe now.

Let $\gamma\in\Omega$ and $\xi$ a vector field along $\gamma$; then
$$
d_\gamma\varphi_h(\xi)=\int_0^1\sigma(\xi(t),\dot\gamma(t))-d_{\gamma(0)}h_t(\xi(t)\,dt=
\int_0^1\langle -J\dot\gamma(t)-\nabla_{\gamma(t)}h_t,\xi(t)\rangle\,dt.
$$
To get this formula, it is sufficient to consider a very narrow film obtained by the translation of $\gamma$ in the direction of $\xi$ and apply the Stokes formula to the integral of $s$ over the boundary of the film.

Hence $\gamma$ is a critical point of $\varphi_h$ if and only if $\dot\gamma=J\nabla_\gamma h_t,\ 0\le t\le 1$.
Note that $J\nabla h_t=\vec h_t$, where $\vec h_t$ is the Hamiltonian vector field associated to the function $h_t$ on the symplectic manifold $M$, i.\,e. $\sigma(\cdot,\vec h_t)=dh_t(\cdot)$. We obtain that critical points of $\varphi_{h_t}$
are exactly periodic solutions of period 1 of the time-varying Hamiltonian system
$$
\dot\gamma(t)=\vec h_t(\gamma(t)),\quad 0\le t\le 1.  \eqno (4)
$$

A similar trick allows to easily compute the Hessian of $\varphi_h$ at a critical point $\gamma$. We take one more vector field $\eta$ along $\gamma$, move our narrow film a little bit in the direction of $\eta$, and apply the Stocks formula to the integral of $\sigma$ over the boundary of the obtained thin 3-dimensional body. We have:
$$
d_\gamma^2\varphi_h(\xi,\eta)=\int_0^1\sigma(\eta(t),\dot\xi(t))-d^2_{\gamma(t)}h_t(\eta(t),\xi(t))\,dt.
$$
It follows that the field $\xi$ belongs to the kernel of $d^2_\gamma\varphi_h$ if and only if $\xi$ is a periodic solution
with period 1 of of the linearization of system (4) along $\gamma$.

Let $P^t:M\to M$ be the flow on $M$ generated by system (4),\,
$$
\frac{\partial P^t}{\partial t}=\vec h_t(P^t(q)),\quad P^0(q)=q,\quad q\in M.
$$
A curve $\gamma(t)=P^t(\gamma(0)),\ 0\le t\le 1$, is a critical point of $\varphi_h$ if and only if $\gamma(0)$ is a fixed point of $P^1$. The linearization of $P^1$ at $\gamma(0)$ is a linear operator ${P^1_*}_{\gamma(0)}: T_{\gamma(0)}M\to T_{\gamma(0)}M$.

We say that $\gamma$ is a non-degenerate 1-periodic solution of system (4) if $\det({P^1_*}_{\gamma(0)}-I)\ne 0$. According to the implicit function theorem, non-degenerate 1-periodic solutions are isolated. Moreover, as we have seen, a 1-periodic solution $\gamma$ of system (4) (a critical point of $\varphi_h$) is non-degenerate if and only if
$\ker D^2_\gamma\varphi_h=0$.

The next step is to associate an index $i(\gamma)\in\mathbb Z$ to any non-degenerate 1-periodic solution of system (4).
To do that we need some general properties of symmetric bilinear forms. Let $b$ be a continuous symmetric bilinear form on a Hilbert space $H$ and $F\subset H$ be a closed subspace. We set
$$
F^{\perp_b}=\{x\in H: Q(x,F)=0\},
$$
the orthogonal complement to $F$ with respect to the form $b$; then $F\cap F^{\perp_b}=\ker(b|_F)$.

Let $F\subset H$ be a finite codimension subspace such that $b|_F$ is nondegenerate\footnote{We say that a bilinear form
$b:H\to H^*$ is non-degenerate if $b$ is invertible.}, then $F^{\perp_b}$ is finite dimensional and the signature
${sgn}(b\bigr|_{F^{\perp_b}})\in\mathbb Z$ is well-defined. Recall that the signature of a symmetric bilinear form is the difference between the positive and negative inertia indices.

\begin{lemma} Let $b_0,b_1$ be symmetric bilinear forms on $H$, $b_s=sb_1+(1-s)b_0$, and $F,\hat F\subset H$, finite codimension subspaces such that $b_s|_F$ and $b_s|_{\hat F}$ are non-degenerate for any $s\in[0,1]$. Then
$$
sgn(b_1\bigr|_{F^{\perp_{b_1}}})-sgn(b_0\bigr|_{F^{\perp_{b_0}}})=
sgn(b_1\bigr|_{\hat F^{\perp_{b_1}}})-sgn(b_0\bigr|_{\hat F^{\perp_{b_0}}}).
$$
\end{lemma}

{\bf Proof. Step 1.} Let $F_0\subset F$ be a subspace of finite codimension such that $b_s|_{F_0}$ is non-degenerate for any $s$ from a subsegment $[\alpha,\beta]\subset[0,1]$. Then
$$
sgn(b_\beta\bigr|_{F_0^{\perp_{b_\beta}}})-sgn(b_\alpha\bigr|_{F_0^{\perp_{b_\alpha}}})=
sgn(b_\beta\bigr|_{F^{\perp_{b_\beta}}})-sgn(b_\alpha\bigr|_{F^{\perp_{b_\alpha}}}).
$$
Indeed, $F_0^{\perp_{b_s}}=F_0^{\perp_{b_s}}\cap F\oplus F^{\perp_{b_s}}$ and $b_s\bigr|_{F_0^{\perp_{b_s}}\cap F}$ is a
non-degenerate form for any $s\in[\alpha,\beta]$. Hence $sgn(b_s\bigr|_{F_0^{\perp_{b_s}}\cap F})$ does not depend on $s\in[\alpha,\beta]$. On the other hand,
$$
sgn(b_s\bigr|_{F_0^{\perp_{b_s}}})=sgn(b_s\bigr|_{F_0^{\perp_{b_s}}\cap F})+
sgn\left(b_s\bigr|_{F^{\perp_{b_s}}}\right).
$$

{\bf Step 2.} The subspace $F\cap\hat F$ has a finite codimension in $H$. Moreover, $\ker(b_s|_{F\cap\hat F})$ are finite dimensional, $0\le s\le 1$. Let $F_s\subset F\cap\hat F$ be a compliment to $\ker(b_s|_{F\cap\hat F})$ in $F\cap\hat F$;
then $b_{s'}|_{F_s}$ is non degenerate for all $s'\in O_s$, where $O_s$ is a neighborhood of $s$ in $[0,1]$. We can cover $[0,1]$ by a finite number of such neighborhoods and apply  Step~1. \qquad$\square$

Lemma~3 allows us to define
$$
\lfloor sgn(b_1)-sgn(b_0)\rfloor\doteq sgn(b_1\bigr|_{F^{\perp_{b_1}}})-sgn(b_0\bigr|_{F^{\perp_{b_0}}});
$$
the defined quantity depends only on $b_1$ and $b_0$ and not on the choice of the subspace $F$. Note that only the difference of the signatures is well-defined; each of the forms $b_0,b_1$ may have infinite inertia indices.

Now consider symmetric bilinear form $d^2_\gamma\varphi_h$. The variables for this form are vector fields along $\gamma$.
The bundle $TM|_\gamma$ is trivializable since $\gamma$ is contractible. Hence one can find an orthonormal basis of
$V_j,W_j, 1\le j \le\frac{\dim M}2$ of $TM|_\gamma$, where $J$ has a canonical form: $JV_j=W_j,\ JW_j=-V_j$. Let $\dim M=2n$; written in the good
basis, $d^2\varphi_h$ turns into a symmetric bilinear form on $H^1(S^1;\mathbb C^n)$ and it has the following expression:
$$
d^2_\gamma\varphi_h(\xi,\eta)=\int_0^1\langle i\eta(t),\dot\xi(t)\rangle-\langle R_t\xi(t),\eta(t)\rangle\,dt. \eqno (5)
$$
Here  $n=\frac{\dim M}2$ and $R_t$ is a real symmetric operator on $\mathbb C^n$ with the $2n\times 2n$-matrix
$\left\{\begin{smallmatrix}d^2_{\gamma(t)}h_t(V_j,V_k) & d^2_{\gamma(t)}h_t(V_j,W_k)\\ d^2_{\gamma(t)}h_t(W_j,V_k) &
d^2_{\gamma(t)}h_t(W_j,W_j)\end{smallmatrix}\right\}_{j,k=1}^n$.

Recall that $\frac d{dt}$ sends the space $H^s(S^1;\mathbb C^n)$ in $H^{s-1}(S^1;\mathbb C^n)$, and
$\int_0^1\langle\cdot,\cdot\rangle\,dt$ is a non-degenerate pairing of $H^s(S^1;\mathbb C^n)$ and $H^{-s}(S^1;\mathbb C^n)$.
It follows that form (5) can be continuously extended to the space $H^{\frac 12}(S^1;\mathbb C^n)$.

From now on we assume that $\xi,\eta\in H^{\frac 12}(S^1;\mathbb C^n)$. Note that
$$
d^2_\gamma\varphi_0-d^2_\gamma\varphi_h=\int_0^1\langle R_t\cdot,\cdot\rangle\,dt
$$
is a compact form on $H^{\frac 12}(S^1;\mathbb C^n)$. Now we analyse the form $d^2_\gamma\varphi_0$

Write the Fourier expansions:
$$
\xi(t)=\sum_{k=-\infty}^{+\infty}e^{i2\pi kt}\xi_k,\qquad \eta(t)=\sum_{k=-\infty}^{+\infty}e^{i2\pi kt}\eta_k.
$$
We set:\quad $\xi_+(t)=\sum\limits_{k=1}^{+\infty}e^{i2\pi kt}\xi_k,\quad \xi_-(t)=\sum\limits_{k=-\infty}^{-1}e^{i2\pi kt}\xi_k$\quad
and similarly for $\eta_+,\, \eta_-$. Then
$$
d^2_\gamma\varphi_0(\xi,\eta)=\sum_{k=1}^{+\infty} 2\pi k(\langle\eta_k,\bar\xi_k\rangle-\langle\eta_{-k},\bar\xi_{-k}\rangle)=\langle\xi_+,\eta_+\rangle_{\frac 12}-
\langle\xi_-,\eta_-\rangle_{\frac 12},
$$
where $\langle\cdot,\cdot\rangle_{\frac 12}$ is the Hilbert inner product of $H^{\frac 12}(S^1;\mathbb C^n)$ (the real part of the Hermitian product).

We have: $H^{\frac 12}(S^1;\mathbb C^n)=H_+\oplus\mathbb C^n\oplus H_-$, where
$H_\pm=\{\xi_\pm: \xi\in H^{\frac 12}(S^1;\mathbb C^n)\}$. Let $s\in[0,1]$; the restriction of the form $d^2_\gamma\varphi_{sh}$ to
$H_\pm$ can be written as follows:
$$
d^2_\gamma\varphi_{sh}(\xi_\pm,\eta_pm)=\pm\langle\xi_\pm+sA_\pm\xi_\pm,\eta_\pm\rangle,
$$
where $A_\pm:H_\pm\to H_\pm$ is a compact symmetric operator

Let $F_\pm\subset H_\pm$ be the orthogonal complement to the linear hull of all eigenvectors of $A_\pm$ corresponding to the eigenvalues that are smaller or equal than $-1$. Then $F=F_+\oplus F_-$ is a finite codimension subspace of
$H^{\frac 12}(S^1;\mathbb C^n)$ and the restriction of the form $d^2_\gamma\varphi_{sh}$ to $F$ is non-degenerate. So we are in the conditions of Lemma~3 and we set:
$$
i_h(\gamma)=\frac 12\lfloor sgn(d^2_\gamma\varphi_0)-sgn(d^2_\gamma\varphi_h)\rfloor. \eqno (6)
$$

Let $\beta_j(M)$ be the Betti number of $M$ of the dimension $j$ (rank of the $j$'s homology group of $M$ with coefficients in a preliminary chosen field) and $C_h$ be the set of all periodic trajectories of system (4) with period 1. If all 1-periodic trajectories are non-degenerate, then $C_h$ is a finite set.

\begin{theorem}[Morse inequalities] Assume that all 1-periodic trajectories of system (4) are non-degenerate. Then, for any $k\in\mathbb Z$, the following inequality holds:
$$
\sum\limits_{j\le k}(-1)^{k-j}\beta_j(M)\le\sum\limits_{\{\gamma\in C_h:i_h(\gamma)\le k\}}(-1)^{k-i_h(\gamma)}.
$$
\end{theorem}

We have $d^2_{\gamma(t)}h(\eta,\xi)=\langle\eta,R_t\xi\rangle,\ \eta,\xi\in T_{\gamma(t)}M$. The linearization of system (4) has a form $\dot\xi=JR_t\xi$. Let $t\mapsto\xi(t)$ be a solution to this equation, then $\xi(t)={P^t_*}_{\gamma(0)}\xi(0),\ 0\le t\le 1$. In what follows, we assume that vector bundle $TM|_\gamma$ is trivialized and treat this differential equation as a linear Hamiltonian system in a fixed symplectic space $T_{\gamma(0)}M$. We also write $T_\gamma M$ and ${P^t_*}_\gamma$ instead of $T_{\gamma(0)}M$  and ${P^t_*}_{\gamma(0)}$ in order to simplify notations a little bit.
Our next goal is to give an effective formula for index (6) in terms of the linear symplectic transformations
${P^t_*}_\gamma:T_\gamma M\to T_\gamma M$. We do it under a mild regularity assumption.

Let $\Sigma,\sigma$ be a symplectic space, $\Pi\subset\Sigma$ a Lagrange subspace, and $L(\Sigma)$ the Lagrange Grassmannian
(the manifold of all Lagrange subspaces of $\Sigma$). Then
$$
\mathcal M_\Pi=\{\Lambda\in L(\Sigma): \Lambda\cap\Pi\ne 0\}
$$
is a codimension 1 cycle in $L(\Sigma)$. Let $\Lambda(t)\in L(\Sigma),\ t_0\le t\le t_1,$ be a curve in $L(\Sigma)$ such that $\Lambda(t_0)\cap\Pi=\Lambda(t_0)\cap\Pi=0$. Recall that the Maslov $\mu_\Pi(\Lambda(\cdot))$
of the curve $\Lambda(\cdot)$ with respect to $\Pi$ is the intersection number of $\Lambda(\cdot)$ and $\mathcal M_\Pi$, see \cite{LiVe} for details.

Now consider symplectic space $\Sigma\times\Sigma$ andowed with the symplectic structure $(-\sigma)\oplus\sigma$.
Given a linear symplectic transformation $Q:\Sigma\to\Sigma$, its graph $\Gamma_Q=\{(\xi,Q\xi):\xi\in\Sigma\}$ is a
Lagrange subspace of $\Sigma\times\Sigma$. In particular, $\Gamma_I$ is the diagonal of $\Sigma\times\Sigma$. Let
$Q(t),\ t_0\le t\le t_1,$ be a curve in the symplectic group such that the maps $Q(t_0)-I$ and $Q(t_1)-I$ are non-degenerate. The Maslov index of $Q(\cdot)$ is just the Maslov index of the curve $\Gamma_{Q(\cdot)}$ with respect to the diagonal,
$$
\mu(Q(\cdot))\doteq\mu_{\Gamma_I}(\Gamma_{Q(\cdot)}).
$$

\begin{theorem} Assume that $\det\left(\int_0^1R_t\,dt\right)\ne 0$. Take $\varepsilon>0$ and consider a curve
$\Phi_\varepsilon$ in the symplectic group defined on the segment $[-\varepsilon,1]$ according to the rule:
$$
\Phi_\varepsilon(t)=\begin{cases} e^{tJ\int_0^1R_\tau\,d\tau},&\text{if $-\varepsilon\le t<0$;}\\
{P^t_*}_\gamma,&\text{if $0\le t\le 1$.}
\end{cases}
$$
Then  $i_h(\gamma)=\mu(\Phi_\varepsilon)$, for any small enough $\varepsilon$.
\end{theorem}

We can make this formula even more explicit assuming more regularity. Indeed, $\mu$ is an intersection number and we can simply count the intersection points. In particular, we assume that $R_t$ continuously depends on $t\in[0,1]$ and introduce subspaces $K_t=\{\xi\in T_\gamma M: R_t\xi=\xi\}\subset T_\gamma M.$
We also denote by $\bar R_t$ the quadratic form on $T_\gamma M$ defined by the formula:
$\bar R(\xi)=\langle R_t\xi,\xi\rangle,\ \xi\in\ T_\gamma M$.

Main regularity assumption is as follows: for any $t\in[0,1]$, if $K_t\ne 0$, then $\bar R_t\bigr|_{K_t}$ is a non-degenerate quadratic form. I emphasize that we require non-degeneracy of the restriction of the form $\bar R_t$
to the subspace $K_t$ and do not care about properties of $\bar R_t$ on whole space $T_\gamma M$. In principle, the form $\bar R_t$ could be non-degenerate on $T_\gamma M$ and be identical 0 on $K_t$. Under this assumption, the points $t\in[0,1]$ such that $K_t\ne 0$ are isolated.

\begin{corollary} Under just imposed assumptions, the following equality holds:
$$
i_h(\gamma)=\sum\limits_{0<t<1}sgn\left(\bar R_t\bigr|_{K_t}\right)+\frac 12\Bigl(sgn(R_0)+sgn\bigl(\int_0^1R_t\,dt\bigr)\Bigr).
$$
\end{corollary}

\section{Sketch of Proofs}

The proof of Theorem~1 is divided in two parts. First we show that for any pair of time-varying Hamiltonians $h_t^0,h_t^1$
and any big enough $c>0$ (how big, depends on our Hamiltonians) there exist $r(c)>0$ and $E(c,r)\in\mathcal E$ such that for any $r>r(c)$ and $E\supset E(c,r)$ the groups $\jmath_r\left(G_{i+d_E}(E;c,r)\right)$ corresponding to $h^0$ and $h^1$
are naturally isomorphic. At this point we do not even need the subspace $E$ to be well-balanced.

In the second part of the proof, we show that $\jmath_r\left(G_{i+d_E}(E;c,r)\right)=H_i(M)$, if $h=0$. Here it is essential that $E$ is well-balanced.

So we start from the choice of $c$. Let us consider the family of time-varying Hamiltonians
$h^\tau_t=\tau h^1_t+(1-\tau)h^0_t,\quad \tau,t\in[0,1].$
The correspondent Hamiltonian fields $\vec h^\tau_t$ are uniformly bounded, hence there exists $c>0$ such that all critical points of $\varphi_{h^\tau}$ belong to the interior of $\varphi^{-1}_{h^\tau}([-c,c]),\ 0\le\tau\le t.$

If $\Omega$ would be a compact finite-dimensional manifold, then we could easily demonstrate homotopy equivalence of the pairs $\bigl(\Omega^c_{h_0},\Omega^{-c}_{h_0}\bigr)$ and $\bigl(\Omega^c_{h_1},\Omega^{-c}_{h_1}\bigr)$. Indeed, the flow
on $\Omega$ generated by a vector field $f_\tau$ such that
$$
\frac{\partial\varphi_{h^\tau}}{\partial\tau}+d_\gamma\varphi_{h^\tau}(f_\tau(\gamma))\le 0,\quad \forall\gamma\in\varphi^{-1}(\pm c),\ 0\le\tau\le 1, \eqno (1)
$$
provides the desired homotopy equivalence.

Unfortunately, $\Omega$ is noncompact and infinite-dimensional. Still $\Omega$ and the differential of $\varphi_h$ are
good enough to guarantee the existence of a finite-dimensional substitute that provides the required weaker equivalence.
I am going to explain how it works.

First of all, $\frac{\partial\varphi_{h^\tau}}{\partial\tau}=\varphi_{h^1}-\varphi_{h^0}$ is uniformly bounded.
Moreover, the gradient of $\varphi_h$ in the $L^2$-norm has a form: $\nabla_\gamma\varphi_h=-J\dot\gamma-\nabla_\gamma h$.
Since any bound on the $L^2$-norm of $\dot\gamma$ implies a bound on $|\varphi_h|$, we obtain that $\nabla_\gamma\varphi_h$
is separated from 0 on $\varphi_h^{-1}(\pm c)$ if $c$ is sufficiently big. Hence we can find an open set of uniformly bounded in the $L^2$-norm by a constant $\rho(c)$ vector fields that satisfy inequality (1).

Now we take $r>0,\ \bar r=r+4\rho(c)$ and a finite-dimensional space $E^l_{\bar r}$ guaranteed by Lemma~1. The weak continuity of the differential of $\varphi_h\circ\phi$ and weak compactness of $B_{\bar r}$ imply the existence of a
finite-dimensional space $E^l\supset E^l_{\bar r}$ and a field $f_\tau$ on $M\times B_{\bar r}$, which takes values in $TM\times E^l$, has the  $L^2$-norm bounded by $2\rho(c)$, and such that
$$
\frac{\partial(\varphi_{h^\tau}\circ\phi)}{\partial\tau}+d_{(q,u)}(\varphi_{h^\tau}\circ\phi)(f_\tau(q,u))<0,
$$
$$
\forall(q,u)\in U_{\bar r}(E)\cap(\varphi_{h^\tau}\circ\phi)^{-1}(\pm c),\quad 0\le\tau\le 1.
$$

The field $f_\tau$ is tangent to the finite-dimensional submanifold $U_{\bar r}(E)$. Moreover, starting in
$U_{r}(E)$ trajectories of this vector field stay inside $U_{\bar r}(E)$ and thus well-defined on the whole
segment $[0,1]$. Moreover, these trajectories stay far enough from the border of $U_{\bar r}(E)$ that allows us to keep necessary control of the homology homomorphisms induced by the imbedding $U_r(E)\subset E\cap\Omega$.

Moreover, $\phi^{-1}(\gamma)$ is an affine subspace of $\gamma(0)\times L^2([0,1];\mathbb C^l)$ for any $\gamma\in\Omega$.
Hence $\phi^{-1}(\gamma)\cap U_{\bar r}(E)$ and $\phi^{-1}(\gamma)\cap U_r(E)$ are convex subsets of $\gamma(0)\times E^l$ (actually, balls) and  $\phi$ induces the homotopy equivalence of the spaces
$U_r(E)\cap\phi^{-1}(\Omega^{\pm c}_{h^\tau})$ and $U_{\bar r}(E)\cap\phi^{-1}(\Omega^{\pm c}_{h^\tau})$
on there images. This is actually all we need to complete the first part of the proof of Theorem~1.

\medskip
Now consider the case $h=0$. Critical points of $\varphi_0$ are exactly constant curves $\gamma(t)\equiv q,\ q\in M$.
In other words, critical points form a smooth manifold $M$. The Hessian of $\varphi_0$ at $q$ is just the quadratic form $b_q$ from
Section~2. As we know, the restriction of this form to a complement to the tangent space to $M$ has zero kernel. If $\Omega$
would be finite-dimensional and $\varphi_0^{-1}([-c,c])$ compact we could say that $\varphi_0$ is a ``Morse--Bott'' function
and conclude that $H_k(\Omega^c_0,\Omega^{-c}_0)=H_{k-\imath}(M)$, where $\imath$ is the negative inertia index of $b_q$,
see \cite{Mi}.

Unfortunately, our objects are infinite-dimensional and noncompact. \linebreak Moreover, both positive and negative inertia indices of
$b_q$ are infinite. Still, as in the first part of the proof, an appropriate finite-dimensional substitute provides all we need.

Given $r>0$ and a sufficiently big well-balanced spaced $E$, let us restrict $\varphi_0$ to $\phi(U_r(E))$ and first study
this restriction in a neighborhood of the manifold $M$ of constant curves. Let $\psi=\varphi_0\bigr|_{\phi(U_r(E))}$ and
$q\in M$; the Hessian of $\psi$ at $q$ is nondegenerate on the complement to the tangent space to $M$ and its positive end negative inertia indices are equal. In particular, both inertia indices are equal to $\frac 12\bigl(\dim(\phi(U_r(E))-\dim M\bigr)$
and there exist a neighborhood of $M$ in $\phi(U_r(E^l))$ that does not contain other critical points of $\psi$.

We need and we have a stronger property: there exists a neighborhood $\mathcal O$ of $M$ in $\Omega$ such that any finite-dimensional subspace of $L^2([0,1];\mathbb C)$ is contained in a well-balanced subspace $E$ such that the restriction
of $\phi_0$ to $\phi(U_r(E))\cap\mathcal O)$ has only constant critical points\footnote{Here ``constants'' means ``constant curves''.}. Indeed, $\nabla_\gamma\varphi_0=-J\dot\gamma$ almost linearly depends on $\dot\gamma$ for $\gamma$ close to the constants and the estimates used to justify Lemma~2 (see, in particular, formula (2.2)) give also a low bound
$|\nabla_\gamma\psi|\ge \epsilon\|\dot\gamma\|$ that is uniform with respect to $E$.

Out of a neighborhood $\mathcal O$, the norm of $\nabla_\gamma\varphi_0$ has a uniform low bound, hence the norm of the gradient-like vector field $\frac{\nabla\varphi_0}{\|\nabla\varphi_0\|^2}$ has an upper bound $\hat\rho$
and we take $r'=r+3c\hat\rho,\ \bar r=r+6c\hat\rho$. Then we find a finite-dimensional approximation $g$ of the field
$-\frac{\nabla(\varphi_0\circ\phi)}{\|\nabla(\varphi_0\circ\phi)\|^2}$ on $(M\times B_{\bar r})\cap\phi^{-1}(\Omega)$ whose norm is bounded by
$\frac 32\hat\rho$ and such that $\langle\nabla(\varphi_0\circ\phi),g\rangle<-\frac 23$, similarly to what we did in the first part of the proof. Finally, we restrict everything to a finite-dimensional manifold $U_{\bar r}(E)$; then any starting in
$U_{r'}(E))\cap\phi^{-1}(\Omega^c_0)$ trajectory of the field $g$ reaches either $\phi^{-1}(\mathcal O)$ or $\phi^{-1}(\Omega^{-c}_0)$, and this is essentially all we need to control homology groups $\jmath_r\bigl(G_i(E;c,r)\bigr)$.

\medskip
Now turn to Theorem~2. We have some freedom in the choice of the generators $X_1,\ldots,X_l$ of $\mbox{Vec}M$ and we are going to adapt them to the periodic trajectories $\gamma\in C_h$. Namely, we may assume that first $\dim M$ generators
give a basis of the bundle $TM|_\gamma$ adapted to the quasi-complex structure $J|_\gamma$. In other words, we may assume
that $d^2_\gamma\varphi_h$ has a form (2.5), where only the ``compact part'' $R_t$ depends on $\gamma\in C_h,\ R_t=R_t^\gamma$.

The proof proceeds similarly to the proof of the second part of Theorem~1. We take $r>0$ big enough and assume that the well-balanced space $E$ is such that $\phi(U_r(E^l))$ contains $C_h$ and all kernels of the bilinear forms
$$
d^2_\gamma\varphi_h(\xi,\eta)=\int_0^1\langle i\eta(t),\dot\xi(t)\rangle-s\langle R^\gamma_t\xi(t),\eta(t)\rangle\,dt,
\quad 0\le s\le 1,\ \gamma\in C_h. \eqno (2)
$$
Moreover, we take $E$ sufficiently big and well-balanced to guarantee that $\varphi_h\bigr|_{\phi(U_{\bar r}(E^l))}$ does not have critical points out of $C_h$ and that the Hessian of $\varphi_h\bigr|_{\phi(U_{\bar r}(E^l))}$ at any $\gamma\in C_h$ is non-degenerate. Then negative inertia index of the Hessian equals
$\frac 12\dim\phi(U_r(E^l))+\imath_h(\gamma)$ and homology groupes $\jmath_r^{\bar r}\bigl(G_i(E;c,r)\bigr)$ are controlled
by these indices and a gradient-like vector field.

\medskip
Theorem~3. To compute the index $\imath_h(\gamma)$ we have to count, according the multiplicities and signs, the
passing through zero eigenvalues of the family of quadratic forms  (2). A field $\xi\in T_\gamma M$ belongs to the kernel
of the form (2) if and only if it is a 1-periodic solution of the linear Hamiltonian system $\dot\xi=sJR_t^\gamma\xi$.
We consider linear symplectic transformations
$$
Q_t(s):\xi(0)\mapsto\xi(t),\quad \xi(0)\in T_\gamma M,
$$
where $\dot\xi(t)=sJR_t^\gamma\xi(t),\ 0\le t\le 1,$ and set $Q(s)\doteq Q_1(s)$.

We see that $\xi$ belongs to the kernel of the form (2) if and only if $Q(s)\xi(0)=\xi(0)$; in other words, if and only if
$\xi(0)\in\Gamma_{Q(s)}\cap\Gamma_I$ (notations of Section~2, just before Theorem~3). The contribution of $Q(s)$ to the
intersection number $\mu(Q(\cdot))$ is equal to the one-half of the signature of the quadratic form
$$
\xi(0)\mapsto\Bigl\langle J\xi(0),\frac{\partial Q}{\partial s}\xi(0)\Bigr\rangle,\quad \xi(0)\in\{\xi\in T_\gamma M: \xi=Q(s)\xi\},
$$
if this quadratic form is non-degenerate. Let $\xi(t)=Q_t(s)\xi(0)$ We have:
$$
\Bigl\langle J\xi(0),\frac{\partial Q(s)}{\partial s}\xi(0)\Bigr\rangle=\Bigl\langle J\xi(1),\frac{\partial Q(s)}{\partial s}\xi(0)\Bigr\rangle=
$$
$$
\int_0^1\Bigl\langle J\dot\xi(t),\frac{\partial Q_t(s)}{\partial s}\xi(0)\Bigr\rangle+\Bigl\langle J\xi(t),\frac{\partial^2 Q_t(s)}{\partial t\partial s}\xi(0)\Bigr\rangle\,dt=
$$
$$
\int_0^1-\Bigl\langle R_t^\gamma\xi(t),\frac{\partial Q_t(s)}{\partial s}\xi(0)\Bigr\rangle+\Bigl\langle\xi(t),R_t^\gamma\frac{\partial Q_t(s)}{\partial s}\xi(0)\Bigr\rangle+\langle\xi(t),R_t^\gamma\xi(t)\rangle\,dt
$$
$$
=\int_0^1\langle\xi(t),R_t^\gamma\xi(t)\rangle\,dt.
$$

On the other hand, according to the standard perturbations theory, the derivative with respect to $s$ of the passing through zero eigenvalue of the form (2) is equal to
$$
\frac\partial{\partial s}\int_0^1\langle i\xi(t),\dot\xi(t)\rangle-s\langle R^\gamma_t\xi(t),\xi(t)\rangle\,dt=
-\int_0^1\langle\xi(t),R_t^\gamma\xi(t)\rangle\,dt.
$$
The formulas for $\imath_h(\gamma)$ and $\mu(Q(\cdot))$ match! We actually got these formulas under a transversality assumption but we can drop this assumption due to the homotopy invariance of the intersection number. We also have:
$Q_0(s)=Q_t(0)=I,\ s,t\in[0,1]$; hence the curve $s\mapsto Qs)=Q_1(s)$ in the symplectic group is homotopic to the curve
$t\mapsto Q_t(1)={P^t_*}_\gamma$ that is presented in Theorem~3, and the Maslov index is homotopy invariant. Additional
``boundary'' terms appear because both curves treated as curves in the Lagrange Grassmannian of the symplectic space $\Sigma\times\Sigma$ start
at the diagonal and the intersection number has to be adjusted according to that.

\section{Step Two Carnot Groups}

In this section we describe an example of an asymptotics of Betti numbers that is much more interesting than in the Floer case and is explicitly computed. This calculation was inspired by the joint work with A.~Lerario and A.~Gentile \cite{AgGeLe},
see also \cite{AgLe,Le}. A motivation, applications and proofs will appear in the forthcoming paper.

A step two Carnot Lie algebra $\frak g$ is a graduated nilpotent Lie algebra with two levels generated by the first level and with a fixed Euclidean structure on the first level. The correspondent simply connected Lie group $\frak G=e^{\frak g}$
is called a step two Carnot group. We have: $\frak g=V\oplus W,\ [V,V]=W,\ [\frak g,W]=0.$
The Euclidean inner product on $V$ is denoted by $\langle\cdot,\cdot\rangle$ and the norm by $|\cdot|$. To any $\omega\in W^*$ we associate an operator $A_\omega\in\mathrm{so}(V)$ by the formula:
$$
\langle A_\omega\xi,\eta\rangle=\langle\omega,[\xi,\eta]\rangle, \quad \xi,\eta\in V. \eqno (1)
$$
It is easy to see that $\omega\mapsto A_\omega,\ \omega\in W^*$ is an injective linear map. Moreover, any injective linear map from $W^*$ to $\mathrm{so}(V)$ defines a structure of step two Carnot Lie algebra on the space
$V\oplus W$ by the same formula read in the opposite direction. We see that step two Carnot Lie algebras are
in the one-to-one correspondence with linear systems of anti-symmetric operators.

Any subspace of Lie algebra forms a left-invariant vector distribution on the Lie group. We are interested in the distribution on the group $\frak G$ formed by $V$. A class $H^1$ curve $\gamma:[0,1]\to \frak G$ is called
{\it horizontal} if $\dot\gamma()\in V_{\gamma(t)}$ for a.\,e. $t\in[0,1]$. Any two points of $\frak G$ can be connected
by a horizontal curve; it follows from the fact that $V$ generates Lie algebra $\frak g$.

Minimum of the lengths of horizontal curves connecting two given points in $\frak G$ is the {\it sub-Riemannian} or
{\it Carnot--Karatheodory} distance between the points. Step two Carnot groups are rather simple and symmetric models of
sub-Riemannian spaces that are not Riemannian.

The following multiplication in $V\times W$ gives a simple realization of $\frak G$ with the origin in $V\times W$ as the unit element:
$$
(v_1,w_1)\cdot(v_2,w_2)=\left(v_1+v_2,w_1+w_2+\frac 12[v_1,v_2]\right).
$$
The Lie algebra is realized by left-invariant vector fields:
$$
(v,w)\mapsto \left(v+\xi,w+\eta+\frac 12[v,\xi]\right),\quad (v,w)\in\frak G,\ \xi\oplus\eta\in V\oplus W=\frak g,
$$
such a field belongs to the horizontal distribution if and only if $\eta=0$.

Starting from the origin horizontal curves are determined by their projection to the first level and have a form:
$$
\gamma(t)=\left(\xi(t),\frac 12\int_0^t[\xi(t),\dot\xi(t)]\,dt\right),\quad 0\le t\le 1,
$$
where $\xi(\cdot)\in H^1([0,1];U),\ \xi(0)=0$. The sub-Riemannian length of $\gamma$ is equal to $\int_0^1|\dot\xi(t)|\,dt$.
For the same reason as in the Riemannian geometry, it is convenient to substitute the length by the essentially equivalent
to the length functional, the action $c\int_0^1|\dot\xi(t)|^2\,dt$, where $c$ is a normalizing constant. For us it is
convenient to take $c=\frac 1{4\pi}$ and we set:
$$
\varphi(\xi)=\frac 1{4\pi}\int_0^1|\dot\xi(t)|^2\,dt.
$$

We focus on the horizontal curves corresponding to closed cures $\xi$; they connect the origin with the second level.
Given $w\in W\setminus 0$, let $\Omega_w$ be the space of horizontal curves connecting $(0,0)$ with $(0,w)$; then
$$
\Omega_w=\left\{\xi\in H^1([0,1];V): \xi(0)=\xi(1)=0,\ \frac 12\int_0^1[\xi(t),\dot\xi(t)]\,dt=w\right\}.
$$
For any $s>0$, we set: $\Omega^s_w=\{\xi\in\Omega_w: \varphi(\xi)\le s\}.$
Note that central reflection $\xi\mapsto -\xi$ preserves $\Omega^s_w$.

Let $E\subset H^1([0,1];V)$ be a finite-dimensional subspace and \linebreak
$\bar E=\bigl(E\setminus0\bigr)/\bigl(\xi\sim(-\xi)\bigr)$ its projectivization, $\bar E$ is homotopy equivalent to
$\mathbb{RP}^{\dim E-1}$. We set $E^s_w=\Omega^s_w\cap E$ and denote by $\bar E^s_w$ the image of $E^s_w$ under the factorization $\xi\sim(-\xi)$.

We consider the homology $H_\cdot(\bar E^s_w;\mathbb Z_2)$ and its image in $H_\cdot(\bar E;\mathbb Z_2)$ by the homomorphism induced
by the imbedding $\bar E^s_w\subset \bar E$. We have:
$$
\mathrm{rank}\bigl(H_i(\bar E^s_w;\mathbb Z_2)\bigr)=\beta_i(\bar E^s_w)+\varrho_i(\bar E^s_w),
$$
where $\beta_i(\bar E^s_w)$ is rank of the kernel of the homomorphism from $H_i(\bar E^s_w;\mathbb Z_2)$ to
$H_i(\bar E;\mathbb Z_2)$ induced by the imbedding $\bar E^s_w\subset\bar E$ and $\varrho_i(\bar E^s_w)$ is the rank of the image of this homomorphism, $\varrho_i(\bar E^s_w)\in\{0,1\}$.

For given $w,E,s$, we introduce two positive atomic measures on the half-line $\mathbb R_+$, the
``Betti distributions'':
$$
\frak b(\bar E^s_w)\doteq \sum\limits_{i\in \mathbb Z_+}\frac 1s\beta_i(\bar E^s_q)\delta_{\frac is},\quad
\frak r(\bar E^s_w)\doteq \sum\limits_{i\in \mathbb Z_+}\frac 1s\varrho_i(\bar E^s_q)\delta_{\frac is}.
$$

Assume that $\dim W=2$ and let $\mathcal E$ be the directed set of all finite-dimensional subspaces of the Hilbert space $H^1([0,1]);V)$.
It appears that there exist limits of these families of measures
$$
\lim\limits_{s\to\infty}\mathcal E\mbox{-}\lim\frak b(\bar E^s_w),\quad
\lim\limits_{s\to\infty}\mathcal E\mbox{-}\lim\frak r(\bar E^s_w)
$$
in the weak topology. Moreover, the limiting measures are
absolutely continuous with explicitly computed densities.

Some notations. Let $\alpha:\Delta\to\mathbb R$ be an absolutely continuous function defined on an interval $\Delta$.
We denote by $|d\alpha|$ a positive measure on $\Delta$ such that $|d\alpha|(S)=\int_S\bigl|\frac{d\alpha}{dt}\bigr|\,dt,\
S\subset\Delta$. Note that the measure $|d\alpha|$ depends only on the function $\alpha$ and not on the choice of the parameter on the interval.

The introduced in (1) operators $A_\omega,\ \omega\in W^*,$ are anti-symmetric and have purely imaginary eigenvalues. Let
$0\le\alpha_1(\omega)\le\cdots\le\alpha_m(\omega)$ are such that $\pm i\alpha_j,\ j=1,\ldots,m,$ are all eigenvalues of
$A_\omega$ counted according the multiplicities; then $\omega\mapsto\alpha_j(\omega)$ are Lipschitz functions.
Let $\bar W^*=(W\setminus 0)/\bigl(w\sim cw,\forall c\ne 0\bigr)$ be the projectivization of $W^*$, $\bar W^*=\mathbb{RP}^1$.

Given $w\in W\setminus 0$, we take the line $w^\perp\in W^*$ and consider the affine line
$$
\ell_w=\bar W^*\setminus\bar w^\perp\subset \bar W^*.
$$
Moreover, we define functions
$$
\lambda^w_j:\ell_w\to\mathbb R_+,\ j=1,\ldots,m,\qquad \phi^w:\ell_w\to\mathbb R_+
$$
by the formulas:
$$
\lambda^w_j(\bar\omega)=\frac{\alpha_j(\omega)}{\langle\omega,w\rangle},\qquad \phi^w(\bar\omega)=\sum_{j=1}^m\lambda^w_j(\bar\omega).
$$
In particular, $\phi^w_*$ transforms measures on $\ell_w$ into measures on $\mathbb R_+$. The Euclidean measure on $\mathbb R_+$ is denoted by $dt$; if $S\subset\mathbb R_+$, then $\chi_Sdt$ is the product of $dt$ and the characteristic function
of $S$.

For simplicity, we compute the limits of the Betti distributions only under a generic assumption on the Carnot group.

\begin{theorem} Assume that there exists $\omega\in W^*$ such that the matrix $A_\omega$ has simple spectrum.
Then, for any $w\in W\setminus 0$, there exist the following limits in the weak topology of the space of positive measures on $\mathbb R_+$:
$$
\frak b_w=\lim\limits_{s\to\infty}\mathcal E\mbox{-}\lim\frak b(\bar E^s_w),\qquad
\frak r_w=\lim\limits_{s\to\infty}\mathcal E\mbox{-}\lim\frak r(\bar E^s_w).
$$
Moreover,
$$
\frak b_w=\phi_*^w\Bigl(\sum\limits_{j=1}^m|d\lambda^w_j|\Bigr), \qquad \frak r_w=\chi_{[0,\min\phi^w]}dt.
$$
\end{theorem}

{\bf Remark.} The internal sequences actually stabilize: there exists $E(s)\in\mathcal E$ such that
$\frak b(\bar E^s_w)=\frak b(\bar E(s)^s_w)$ and $\frak r(\bar E^s_w)=\frak r(\bar E(s)^s_w)$ for any $E\supset E(s)$.
An interesting feature of the Betti measures $\frak b_w$ and $\frak r_w$ is their sensitivity to the endpoint $w$ that reflects well the anisotropy of the space of horizontal curves already at the homological level, even if the spaces $\Omega_w$ are all contractible.

\section{Conclusion}

We conclude with a draft of general scheme for a soft construction of Floer-type asymptotic homologies.
The object to study is a Banach manifold $\Omega$ equipped with a growing family of closed subsets
$\Omega^s,\ s\in\mathbb R$. Auxiliary objects are a Banach space $B$ and a submersion $\Phi:U\to \Omega$, where
$U\subset B$ is a finite codimension submanifold of $B$. Moreover, $U$ is equipped with an ordered by the inclusion directed and exhausting family $\mathcal V$ of open bounded subsets and $B$ is endowed by an ordered by the inclusion directed family
$\mathcal E$ of finite dimensional subspaces such that $\overline{\bigcup\limits_{E\in\mathcal E}E}=B$.

Given $E\in\mathcal E,\ V\in\mathcal V,\ s\in\mathbb R,$ we consider the relative homology groups:
$$
G_i(E,V,s)\doteq H_i\left(\Omega^s\cap\Phi(E\cap V),\Omega^{-s}\cap\Phi(E\cap V)\right).
$$
Moreover, for $V\in\mathcal V$ we denote by
$$
\jmath_V:G_i(E,V,s)\to H_i\left(\Omega^s\cap\Phi(E\cap U),\Omega^{-s}\cap\Phi(E\cap U)\right)
$$
the homology homomorphism induced by the inclusion $V\subset U$.

Finally, we select normalizing quantities
$r_i(E,s),\,\rho_i(E,s)\in\mathbb R_+$ and build atomic measures:
$$
\frak b(E,V,s)=\sum_{i\in\mathbb Z_+}\rho_i(E,s)\mathrm{rank}\bigl(\jmath_VG_i(E,V,s)\bigr)\delta_{r_i(E,s)}
$$
in such a way that their exist a nonzero limit:
$$
\frak b=\lim_{s\to\infty}\mathcal V\mbox{-}\lim\mathcal E\mbox{-}\lim\frak b(E,V,W,s).
$$
We treat this limit as an asymptotic distribution of Betti numbers.

Some of ingredients may be superfluous. In the Leray--Schauder case we take a single limit (variables $s\in\mathbb R$ and
$V\in\mathcal V$ are absent) and in the two step
Carnot group case we take the double one (variables $V\in\mathcal V$ is absent). Of course, the equivariant version is also available if there is a natural group action; in the Carnot group case we used the involution $\xi\mapsto(-\xi)$.

\end{document}